%

\input mssymb


%

%

\newcount\skewfactor
\def\mathunderaccent#1#2 {\let\theaccent#1\skewfactor#2
\mathpalette\putaccentunder}
\def\putaccentunder#1#2{\oalign{$#1#2$\crcr\hidewidth
\vbox to.2ex{\hbox{$#1\skew\skewfactor\theaccent{}$}\vss}\hidewidth}}



\def\rest{\mathord{\restriction}}

\def\phi{\varphi}

\def\su{\subseteq}
\def\a{\alpha}
\def\b{\beta}

\def\l{\lambda}
\def\k{\kappa}

\def\om{\omega}

\def\lng{\langle}
\def\rng{\rangle}

\def\sm{\setminus}
\def\cont{{2^{\aleph_0}}}
\def\limdir{{\lim\limits_{\longrightarrow}}}
\def\liminv{{\lim\limits_{\longleftarrow}}}

\def\discup{\overnext{\bullet}\cup}
\def\bigdiscup{{\bigcup\limits^{\bullet}}}


\def\dom{{\rm dom}}

\def\ran{{\rm  ran}}
\def\equiv{\sim}

\def\id{{\rm id}}

\def\aut{{\rm Aut}}

\def\fin{{\rm Fin}}


\def\endproof#1{\hfill  
{\parfillskip0pt$\smiley_{\hbox{{#1}}}$\par\medbreak}}

\def\iff{\Leftrightarrow}
\def\proof{\smallbreak\noindent{\sl Proof}: }

\def\conc{\concatenate}


\def\Cal#1{{\cal #1}}


\newbox\noforkbox \newdimen\forklinewidth
\forklinewidth=0.3pt   

\setbox0\hbox{$\textstyle\bigcup$}
\setbox1\hbox to \wd0{\hfil\vrule width \forklinewidth depth \dp0
			height \ht0 \hfil}
\wd1=0 cm
\setbox\noforkbox\hbox{\box1\box0\relax}
\def\unionstick{\mathop{\copy\noforkbox}\limits}
\def\nonfork#1#2_#3{#1\unionstick_{\textstyle #3}#2}
\def\nonforkin#1#2_#3^#4{#1\unionstick_{\textstyle #3}^{\textstyle  
#4}#2}

\setbox0\hbox{$\textstyle\bigcup$}
\setbox1\hbox to \wd0{\hfil$\nmid$\hfil}
\setbox2\hbox to \wd0{\hfil\vrule height \ht0 depth \dp0 width
				\forklinewidth\hfil}
\wd1=0cm
\wd2=0cm
\newbox\doesforkbox
\setbox\doesforkbox\hbox{\box1\box0\relax}
\def\nunionstick{\mathop{\copy\doesforkbox}\limits}

\def\fork#1#2_#3{#1\nunionstick_{\textstyle #3}#2}
\def\forkin#1#2_#3^#4{#1\nunionstick_{\textstyle #3}^{\textstyle  
#4}#2}

\font\circle=lcircle10

\setbox0=\hbox{~~~~~}
\setbox1=\hbox to \wd0{\hfill$\scriptstyle\smile$\hfill} 
\setbox2=\hbox to \wd0{\hfill$\cdot\,\cdot$\hfill} 

\setbox3=\hbox to \wd0{\hfill\hskip4.8pt\circle i\hskip-4.8pt\hfill}  


\wd1=0cm
\wd2=0cm
\wd3=0cm
\wd4=0cm

\newbox\smilebox
\setbox\smilebox \hbox {\lower 0.4ex\box1
		 \raise 0.3ex\box2
		 \raise 0.5ex\box3
		\box4
		\box0{}}
\def\smiley{\leavevmode\copy\smilebox}

\headline={\tenrm  
\number\folio\hfill\jobname\hfill\number\day.\number\month.\number\year}


\outer\long\def\ignore#1\endignore{}

\newcount\itemno
\def\itm{\advance\itemno1 \item{(\number\itemno)}}
\def\ritm{\advance\itemno1 \item{)\number\itemno(}}
\def\startitm{\itemno=0 }
\def\aitm{\advance\itemno1 

\item{(\letter\itemno)}}
^^L

\def\letter#1{\ifcase#1 \or a\or b\or c\or d\or e\or f\or g\or h\or
i\or j\or k\or l\or m\or n\or o\or p\or q\or r\or s\or t\or u\or v\or
w\or x\or y\or z\else\toomanyconditions\fi}
\def\raitm{\advance\itemno1 \item{)\rletter\itemno(}}
\def\rletter#1{\ifcase#1\or `\or a\or b\or c\or d\or e\or f\or g\or
h\or i\or k\or k\or l\or n\or n\or q\or r\or t\or v
\else\toomanyconditions\fi}


\newcount\secno
\newcount\theono

\catcode`@=11
\newwrite\mgfile

\openin\mgfile \jobname.mg
\ifeof\mgfile \message{No file \jobname.mg}
	\else\closein\mgfile\relax\input \jobname.mg\fi
\relax
\openout\mgfile=\jobname.mg

\newif\ifproofmode
\proofmodetrue            

\def\@nofirst#1{}

\def\neusection{\advance\secno by 1\relax \theono=0\relax}
\def\neuchap{\secno=0\relax\theono=0\relax}

\neuchap

\def\labelit#1{\global\advance\theono by 1%
             \global\edef#1{%
             \number\secno.\number\theono}%
             \write\mgfile{\@definition{#1}}%
}

^^L




\def\ppro#1#2:{%
\labelit{#1}%
\smallbreak\noindent%
\@markit{#1}%
{\bf\ignorespaces {#2}:}}





\def\@definition#1{\string\def\string#1{#1}
\expandafter\@nofirst\string\%
(\the\pageno)}

\def\@markit#1{
\ifproofmode\llap{{ \expandafter\@nofirst\string#1\ }}\fi%
{\bf #1\ }
}

\def\h@markit#1{
\ifproofmode\edef\nxt{\string#1\ }%
{\tenrm\beginL\nxt\endL}
\fi%
{{\bf\beginL #1\endL}}
}
 ^^L
\def\labelcomment#1{\write\mgfile{\expandafter
		\@nofirst\string\%---#1}}

\catcode`@=12

\def\refrence#1#2:{\write\mgfile{\def\noexpand#1{#2}}%
\areference{#2}}
\def\areference#1{\medskip\item{[#1]} \ignorespaces}

\newcount\referencescount
\def\numrefrence#1{\advance\referencescount1
\edef#1{\number\referencescount}%
\write\mgfile{\def\noexpand#1{#1}}%
\areference{#1}}

\def\numericalreferences{\let\refrence\numrefrence}

\newcount\scratchregister
\def\simplepro{\scratchregister\theono\advance\scratchregister1 

\edef\scratchmacro{\number\secno.\number\scratchregister}%
\expandafter\ppro\csname\scratchmacro\endcsname}


\magnification=1200
\baselineskip=16pt

\def\Sym{{\rm Sym}}

\def \F{{\Cal F}}
\def\G{\Gamma}
\def\FG{{\rm FG}}
\def\emptyseq{\lng\rng}
\def\limdir{{\lim\limits_{\longrightarrow}}}
\def\discup{\dot\cup}
\def\bigdiscup{{\bigcup\limits^{.}}}
\def\bigdis{\bigdiscup}
\def\conc{\circ}
\def\m#1#2{\mathord{\mathop{#2}\limits^{#1}}{}}
\def\pi{\Phi}
\overfullrule0pt
\headline{}

\font\bigfont cmbx10 scaled \magstep2
\font\namefont cmbx10 scaled \magstep2
{{
\obeylines

\everypar={\hskip0cm plus 1 fil}
\bf
\parskip=0.4cm

{\bigfont HOMOGENEOUS FAMILIES}
{\bigfont AND}
{\bigfont THEIR AUTOMORPHISM GROUPS}

\bigskip
{\rm May 1993}
\bigskip

\parskip0.1cm

{\namefont Menachem Kojman}\footnote{$^{*}$}{\rm Partially supported  
by  

the Edmund Landau Center for research in Mathematical Analysis,  
sponsored by
the Minerva Foundation (Germany).}
Department of Mathematics 

Carnegie Mellon University
Pittsburgh, PA 15213, USA
{\tt kojman@andrew.cmu.edu}

\bigskip

{\namefont Saharon Shelah}\footnote{$^{**}$}{\rm The Second author  
thanks the
Israeli Academy of Sciences for partial support. Publication number  
499}
Institute of Mathematics 

Hebrew  University of Jerusalem, Givat Ram
91904 Jerusalem, Israel
Rutgers University, 

New Brunswick NJ, USA
{\tt shelah@math.huji.ac.il}

}
\vfill
\everypar{}
\rm
\leftskip2cm\rightskip2cm

ABSTRACT. A homogeneous family of subsets over a given set is one  
with
a very ``rich'' automorphism group. We prove the existence of a
bi-universal element in the class of homogeneous families over a  
given
infinite set and give an explicit construction of $2^\cont$
isomorphism types of homogeneous families over a countable set.

\footline{\hfill}
\global\pageno0

\eject
}

\proofmodefalse

\ignore 
A bipartite graph $\G$, is a triple $(L,R,E)$ where $L$ is the
set of left vertices, $R$ that of the right vertices, $L\cap
R=\emptyset$ and $E\su [L\cup R]^2$ is the set of edges, and $E\cap
[L]^2=E\cap [R]^2=\emptyset$. We write $xEy$ for $\{x,y\}\in E$. By
an {\it
automorphism} of $\G$ we always mean a 1-1 and onto function from
$L\cup R$ to $L\cup R$ which preserves edges {\it and sides}, that
is $L$ is mapped to $L$ and $R$ to $R$. We call $\G$ {\it
homogeneous} if every finite partial automorphism of $\G$ can be
extended to an automorphism.

A bipartite graph $\G$ is said to be {\it of type $(\k,\l)$} if the
cardinality of its left side is $\k$ and the cardinality of its
right side is $\l$. In this paper we discuss homogeneous bipartite
graphs of type $(\aleph_0,\cont)$. We shall prove:

{\bf A. Theorem}:  There is a {\it bi-universal} homogeneous
bipartite
graph of type $(\aleph_0,\cont)$,  that is one in which all
bipartite graphs of type $(\aleph_0,\cont)$, not necessarily
homogeneous, are embedded as induced subgraphs.

{\bf B. Theorem}: There are $2^\cont$ isomorphism types of
homogeneous bipartite graphs of type $(\aleph_0,\cont)$.

We shall use the method of   direct limits to construct homogeneous
bipartite graphs.

 The existence of homogeneous bipartite graphs of type
$(\aleph_0,\cont)$ was proved in [GGK], where is was also shown
that there are $2^{\aleph_1}$ isomorphism types of homogeneous
graphs of type $(\aleph_0,\cont)$ if the continuum hypothesis
holds, and that there is a  unique isomorphism type of a
homogeneous $(\aleph_0,\aleph_1)$ bipartite graphs if CH fails and
MA holds.

AUTOMORPHISM GROUPS?
TOPOLOGY?


We say that a bipartite graph $\G=(L,R,E)$ is {\it extensional} if
$E$ satisfies the extensionality axiom $(\forall
x,y)[x=y\leftrightarrow (\forall z)(xEz\leftrightarrow
yEz)]$. In other words, an extensional bipartite graph (in
fact any graph) $\G$ is one in which a vertex $x$ is determined by
its {\it set of neighbors} $\G(x)=\{y\in \G:xEy\}$. An
extensional bipartite graph is naturally isomorphic to an {\it
incidence graph} of a family of sets: think of the left side $L$
as { points} and identify each vertex in $R$ with its set of
neighbors. Thus we have specified a family $\F$ of subsets of
$L$. The incidence graph of the family is $(L,\F,\in)$,  and is
clearly isomorphic to $\G$.

If $\G$ is a homogeneous bipartite graph, we argue that it must be
extensional if it is not trivial:

\ppro \extens Lemma: {\sl Suppose that $\G=(L,R,E)$ is a homogeneous
bipartite graph and $E\not=\emptyset, E\not=L\times R$. Then $\G$
is extensional.}

\proof: Suppose that $x,y$ are distinct vertices of the same side,
without loss of generality 

$L$, and $\G(x)=\G(y)$.
Let $z\in L$ be any  vertex, and consider the partial
automorphism that maps $x$ to $x$ and $y$ to $z$. As this partial
automorphism extends to a total automorphism, it follows that
$\G(x)=\G(z)$. Therefore $\G(x)=\G(z)$ for all $z\in L$. If
$\G(x)=R$, then $E=L\times R$, and if $\G(x)=\emptyset$ then
$E=\emptyset$. Otherwise, there is $v\in \G(x)$ and $u\in R\sm
\G(x)$. Clearly, there is no automorphism that carries $u$ to
$v$ (as $u$ has neighbors, while $v$ has none). This
contradicts homogeneity.\endproof\extens

We exclude graphs with either empty or full set of edges from our
discussion, and therefore, by the lemma above, consider only
extensional bipartite graphs. We  shall confine our discussion,
then, to incidence graphs of families of sets. 

\endignore 

\beginsection{\S0 Introduction}

Homogeneous objects are often defined in terms of their automorphism
groups. Rado's graph $\Gamma$, also known as the countable random
graph, has the property that for any isomorphism $f$ between two
finite induced subgraphs of $\Gamma$ there is an automorphism of
$\Gamma$ extending $f$. This property is the homogeneity of Rado's
graph; and any graph whose automorphism group satisfied this  
condition
is called homogeneous.

The automorphism group of Rado's graph was studied by Truss in [T2],
and shown to be simple. Truss studied also the group AAut$(\Gamma)$  
of
{\it almost automorphisms} of Rado's graph (see [T3] and also  
[MSST]).
This is a {\it highly transitive} group extending Aut$(\Gamma)$
(where  ``highly transitive'' stands for  ``$n$-transitive for all  
$n$'';
the group $\aut(\Gamma)$ is not highly transitive).

In this paper we shall study {\it homogeneous families of sets} over
infinite sets. Our definition of homogeneity of a family of sets
implies that its automorphism group satisfies, among other  
conditions,
that it is highly transitive. However, while all homogeneous graphs
over a countable set are classified (see [LW]), this is not the case
with homogeneous families over a countable set.

We shall show that there are $2^{2^{\aleph_0}}$ isomorphism types of
homogeneous families over a countable  set. This is done in Section  
4.
From the proof we shall get $2^{2^{\aleph_0}}$ many permutation
groups, each acting homogeneously on some family over $\om$, and each  
being
isomorphic to the free group on $2^{\aleph_0}$
generators, but such that no two are conjugate in $\Sym(\om)$.

In Section 3 we prove the existence of a bi-universal homogeneous
family over any given infinite set. The definitions of bi-embedding
and bi-universality are generalization of definitions made by Truss  
in
his study of universal permutation groups [T1]. A short survey of
results concerning the existence of universal objects can be found in
the introduction to [KS1]. Results concerning abelian
groups are in [KS2], and results on stable unsuperstable first
order theories are in [KS3].

Homogeneous families were studied in [GGK] (where they were treated  
as
bipartite graphs). There it was shown that the number of isomorphism
types of homogeneous families over $\om$ of size $\aleph_1$ is
independent of ZFC and may be 1 as well as $2^{\aleph_1}$ in  
different
models of set theory.

Model theorists will recognize that uncountable homogeneous families
over a countable set are examples of two-cardinal models which are
$\om$-homogeneous as well. Set theorists may be interested in the
following

\ppro \problem Problem: Is it consistent that $2^{\aleph_0}$ is large
and that in some uncountable $\l<2^{\aleph_0}$ there is a maximal
homogeneous family (with respect to inclusion)?

We wish to remark finally that the existence of $2^\cont$ isomorphism
types of homogeneous families over $\om$ follows from a general
theorem about non-standard logics [Sh-c, VIII,\S1] (for more details
see also [Sh 266]). The virtue of the proof here (besides being
elementary) is its explicitness and the information it gives about  
the
embeddability of an arbitrary family in a homogeneous one.

{\bf NOTATION} We denote disjoint unions by $\discup$ and $\bigdis$.  
A
natural number $n$ is the set $\{0,1,\ldots n-1\}$ of all smaller
natural numbers.

\neusection
\beginsection{\S1 Getting started}

Let $\F\su\Cal P(A)$ be a family of subsets of a
given  infinite set $A$. An {\it automorphism} of $\F$ is a
permutation $\sigma\in
\Sym(A)$ which  satisfies that $X\in
\F\iff \sigma[X]\in \F$ for every $X\su
A$. (By $\sigma[X]$ we denote $\{\sigma(x):x\in X\}$ for $X\su A$.)  
The
group $\aut(\F)\su \Sym(A)$ is the group of all automorphisms of  
$\F$.

One way of defining when a family $\F\su\Cal P(A)$ is homogeneous is  
to
demand that the bipartite graph $\lng A,\F,\in\rng$ is homogeneous,
namely that every finite partial automorphism of this graph which  
respects
the sides extends to a total automorphism. We shall write a more
complicated (though equivalent) definition. This will be needed in
what follows.

\ppro \demand Definition: {\sl Suppose $\F\su \Cal P(A)$ is a given
family of subsets of a set $A$. A {\it demand on $\F$} is a pair
$d=(h^d,f^d)$ such that $h^d$ is a finite 1-1 function from $A$ to
$A$, $f^d$ is a finite 1-1 function from $\F$ to $\F$ and $x\in
X\iff  h^d(x)\in f^d(X)$ for every $x\in \dom h^d, X\in \dom f^d$.
We denote by $D=D(A,\Cal F)$ the set of all demands on $\F$. Let
$\FG(D)$ be the free group over the set $D(A,\F)$.   We say that an
automorphism $g\in \aut(\F)$ {\it satisfies} a demand $d$ if
$g(x)=h^d(x)$ for $x\in \dom h^d$ and $g[X]=f^d(X)$ for $X\in \dom
f^d$. 

\ignore
Let $S\su D$ be the set of all {\it satisfiable} demands,
namely
$S=\{d\in D:(\exists g\in \aut(\F))(g \;\hbox{\ satisfies\ }\; d)\}$.
\endignore
We
call a partial homomorphism $\phi:\FG(D)\to \aut(\F)$ a {\it
satisfying}
homomorphism if $\phi(d)$ satisfies $d$ for $d\in \dom \phi$. (By
``partial'' we mean that $\phi$ need not be defined on all generators
of $\FG(D)$.)}

\ppro \rephrase Definition: {\sl A family $\F$ is {\it homogeneous}  
if
and only if every $d\in D$ is satisfiable if and only if there is a
(total) satisfying homomorphism $\phi:\FG(D)\to \aut(\F)$. A group
$G\su \aut \F$ {\it acts homogeneously} on $\F$ if and only if $G$
contains the image of a total satisfying homomorphism, or,
equivalently if and only if every demand is satisfied by some element
in $G$.}

When $\phi$ is a homomorphism as above, we say that $\phi$ {\it
testifies} the homogeneity of $\F$.

  When the set
$A$ is clear from the context, we write $D(\Cal F)$ instead of
$D(A,\Cal F)$.

\ppro \example Examples:  {\sl
\startitm
\itm The family $\F=\{\{x\}:x\in A\}$ of all
singletons is homogeneous. The group $\aut(\F)$ is the group
$\Sym(A)$ of all symmetries of $A$.
\itm The family $\fin(A)$ of all finite subsets of $A$ is not
homogeneous, although $\aut(\fin(A))=\Sym(A)$, because a demand
$d=(\emptyset, \{(X,Y)\})$ cannot be satisfied when $X$ and $Y$ are
finite sets of different cardinalities. 

\itm A countable family of random subsets of
$\om$ is homogeneous in probability 1. The membership of a point in   
a
random set is determined by flipping a coin.}

In [GGK] the following was proved:

\ppro \classification Theorem:  {\sl Every homogeneous family of  
subsets of
an infinite set $A$ satisfies exactly one of the conditions below:

\startitm
\itm $\F=\{\emptyset\}$
\itm $\F=\{A\}$
\itm $\F$ is the family of all singletons of $A$
\itm $\F$ is the family of all co-singletons of $A$
\itm $\F$ is an independent family, namely for every finite function
$\tau:\F\to \{+,-\}$, the set $B_\tau=\bigcap _{X\in\tau^{-1}(+)}X\
\cap \ \bigcap_{Y\in \tau^{-1}(-)} A\sm Y$ is infinite, and $\F$ is
dually independent, namely for every function $\tau: A\to \{+,-\}$
there are infinitely many members of $\F$ containing $\tau^{-1}(+)$  
and
avoiding $\tau^{-1}(-)$. Equivalently, the first order theory of  
$\lng
A,\F,\in\rng$ is the first order theory of the random countable
bipartite graph.  }

\ignore

We shall explicitly construct  $2^\cont$ many isomorphism
types of (independent) homogeneous families over a countable set $A$  
in Section 3.

Aside from counting the number of non isomorphic families, we are
interested in the embeddability relation among homogeneous families,
and particularly in the existence of a {\it universal } homogeneous
family, that is one in which all others are embeddable.

One natural definition for ``the family $\F_0\su \Cal P(A_0)$ is
embeddable into the family$\F_1\su\Cal P(A_1)$'' is: ``there is a 1-1
onto function $f:A_0\to A_1$ such that $f[X]\in\F_1$ for all $X\in
\F_0$''. When $f$ is identity, this  means that $\F_0$ is a
subfamily of $\F_0$.  Under this definition there is no
universal homogeneous family over a countable set (see Claim \noUni\
below).

Another natural definition for embeddability of families arises from
the usual model theoretic definition of embeddability. View a family
$\F\su \Cal P(A)$ as a bi-partite graph, $\lng A,\F,E\rng$ where
$(x,X)\in E\iff x\in X$. Now say that $\F_0$ is embedded in $\F_1$ if
the bi-partite graph associated with $\F_0$ is isomorphic to an
induced subgraph of the bi-partite graph associated with $\F_1$.

This second definition is equivalent to saying that there is a 1-1
function $f:A_0\to A_1$ which is not necessarily onto, such that for
every $X\in\F_0$ there is $Y\in
\F_1$ such that $f[X]=\{x\in Y\cap \ran f\}$. We will show that there
is a universal homogeneous family with respect to this type of
embedding (theorem \universal\ below) and something more.

We say that a family $\F_0\su\Cal P(A_0)$ is {\it bi-embedded} in
$\F_1\su\Cal P(A_1)$, or ``embedded together with its automorphism
group'' if there is a 1-1 function $f:A_0\to A_0$ which is an
embedding of $\F_0$ in $\F_1$ and such that every automorphism of the
image of $\F_0$ under $f$ is the restriction to the range of $f$ of
some automorphism of $\F_1$. 

When a family $\F_0$ is bi-embedded in
$\F_1$, it is isomorphic to a subfamily of all restrictions $\{Y\cap
A:Y\in \F_1\}$ for a fixed $A\su A_1$, and the group $\aut(\F_0)$ is
isomorphic to a subgroup of the set-wise stabilizer
$\aut(\F_1)_A=\{\sigma\in\aut(\F_1):\sigma\rest A\in\Sym(A)\}$.

The main result we prove in Section 2 is the existence of a  
bi-universal
homogeneous family, namely a homogeneous family in which all others
are bi-embedded. 

\vfill
\eject

\noindent
{\bf Various  Contexts}

We address first a few words to graph theorists. Homogeneous families
were dealt with in [GGK] in the attempt to answer the following
``local/global'' type question about bipartite graphs: Is it true  
that
in every connected bi-partite graph of girth 4 which is locally
3-symmetric both sides are of equal cardinality?  This was proved 

for {\it finite} bipartite graphs by Perles and Kupitz, who asked if
the assumption of finiteness could be removed. Uncountable  
homogeneous
families show that the finiteness assumption is essential even when
strengthening the local symmetry requirements from local 3-symmetry  
to
homogeneity.

Next we wish to say a few words to model theorists. The themes of
counting isomorphism types and studying the relation of embeddability
are common in model theory in general (and in the work of the second
author in particular). An uncountable homogeneous family, the model
theorist will notice, is a two-cardinal model which is also
$\om$-homogeneous.  In the class of such models many model theoretic
methods are not useful (most notably the compactness theorem). Hence
the methods presented here may be of interest for the model theorist.

The existence of $2^\cont$ isomorphism types of homogeneous families
follows from a general theorem about non-standard logics [Sh-c,
VIII,\S1] (for more details see also [Sh 266]). The virtue of the
proof here (besides being elementary) is its explicitness and the
extra information it gives about the embeddability of an arbitrary
family in a homogeneous one.

Finally, we have a word or two for set theorists of the continuum.
Except for trivial cases, all homogeneous families are independent
families. Let $HI$ be the set of cardinalities in which there are  
{\it
maximal} homogeneous independent families.  It is not clear that
$BI\not=\emptyset$. Define $\frak hi$ as the minimal member of $HI$
when $HI$ is not empty.  The continuum hypothesis and also Martin
Axiom imply that ${\frak hi}={\frak c}$ (See [GGK]). We urge the set
theorists of the continuum to absorb $\frak hi$ into the community of
cardinal invariants of the continuum, and in particular to decide
whether it may (exist and) be different from the good old $\frak i$
(the minimal cardinality of a maximal independent family).

We conclude with a few words on universal objects. The question of  
the
existence of universal members in various categories is natural, and
arises independently in many areas of mathematics. A short survey of  
results
concerning the existence of universal objects can be found in the
introduction to [KjSh 409]. Additional results concerning abelian
groups is in [KjSh 455], and results on stable unsuperstable first
order theories are in [Kjsh 447]. 

\medbreak
\noindent
\endignore

\neusection
\beginsection{\S2 Direct limits and homogeneity}

In this section we exhibit a method of constructing homogeneous
families as direct limits. This method will be used in the following
sections.

Homogeneity is not, in general, preserved under usual direct limits  
of
families. For example, an increasing union of homogeneous families
need not be homogeneous itself. We therefore consider here a stronger
relation of embeddability,  called here ``multi-embeddability'',  
which,
roughly speaking, preserves the satisfaction of  previously satisfied
demands. Direct limits of {\it this} relation can be made  
homogeneous, as we
shall presently see.
\medskip

\ppro \embedding Definition: {\sl Let $T_i=\lng A_i,\F_i,
D_i,G_i,\phi_i\rng$, ($i=0,1$), be respectively  a set $A_i$, a
family
of
subsets $\F_i\su\Cal P(A_i)$, the collection of demands
$D_i=D(\F_i)$, an automorphism group  $G_i\su\aut(\F_i)$ and a
partial
satisfying
homomorphism $\phi_i:\FG(D_i)\to G_i$. Let $\bar T_i= A_i\cup
\F_i\cup
D_i\cup G_i$. We
call a function $\pi:\bar T_0\to \bar T_1$  a  {\it multi-embedding
of $T_0$ into $T_1$} (and write $\pi:T_0\to T_1$) if:
\itm $\pi\rest A_0$ is a 1-1 function into $A_1$
\itm $\pi\rest \F_0$ is a 1-1 function into $\F_1$
\itm $\pi\rest D_0$ is a 1-1 function into $D_1$ 

\itm $\pi\rest G_0$ is a group monomorphism into $G_1$

And the following rules hold for $x\in A_0$, $X\in \F_0$, $d\in
D_0$ and $g\in G_0$:
\startitm
\aitm $x\in X\iff \pi(x)\in \pi(X)$
\aitm $\pi[\dom h^d]= \dom h^{\pi(d)}$, $\pi[\dom
f^d]= \dom f^{\pi(d)}$,
$\pi((h^d(x))=h^{\pi(d)}(\pi(x))$ and 

$\pi((f^d(X))=f^{\pi(d)}(\pi(X)$ 

\aitm $\pi(g(x))=\pi(g)(\pi(x))$ and $\pi(g[X])=\pi(g)[\pi(X)]$ 

\aitm  $\pi(d)\in\dom\phi_1$ and  $\pi(\phi_0(d))=\phi_1(\pi(d))$ 

for
every $d\in \dom\phi_0$

We say that a multi-embedding $\pi$ is {\it successful} if in
addition to
the conditions above also the following holds
\aitm $\pi(d)\in\dom\phi_1$ for every $d\in D_0$.}

\ignore 

It is convenient to assume that $pi\rest A_0=\id$. Under this
assumption the conditions above say that $A_0$ is a subset of $A_1$
which is invariant under $\pi[G_0]$, that $X=\pi(X)\cap A_0$ for
every
$X\in \F_0$ and that $g=\pi(g)\rest A_0$ for every $g\in G_0$.

{\bf Discussion} Coming to satisfy the demands in $D_1$ by
automorphisms in $G_2$, there are those demands which ``arrived''
from
$D_0$ and for which an automorphism was already
chosen by $\phi_0^1$. The coherence condition ensures that this
choice
is not
changed, of, more precisely, that $\phi_0^2(d)$ is compatible with
$\phi_0^1(d)$ for $d\in d_0$, namely $\phi_0^2(d)\rest
A_1=\phi_0^1(d)$, and --- even more --- that
$\phi_0^2(d)=\pi_1^2(\phi_0^1(d))$.

We shall be forming  direct limits of successful and coherent
embeddings in the intention that in the limit all demands are
satisfied.  

\endignore 

\ppro \direct Definition:  {\sl
Suppose $I$ is a directed set  and $T_i=\lng
A_i, \F_i, D_i,G_i,\phi_i\rng$ is as in definition \embedding\ above
for
$i\in I$.  Suppose that $\pi_i^j:T_i \to T_j$ is a multi-embedding
for
$i\le j$, and
\item{(i)} $\pi_i^i=\id$ 

\item{(ii)} $\pi_j^k\pi_i^j=\pi_i^k$ for $i\le j\le k$.

Then we call ${\bf T}=\lng T_i: (i\in I);(\pi_i^j,\phi_i^j)\rng$  a
{\it direct system of  multi-embeddings}. We call ${\bf T}$ {\it
successful}
if in addition to (i) and (ii) the following condition holds:

\item{(iii)} for every $i\in I$ there is $j\ge i$ such that $\pi_i^j$
is  successful. }

\ppro \DirIsHom Theorem:  {\sl Suppose $\lng T_i:(i\in
I);\;(\pi_i^j,\phi_i^j)\rng$ is a successful direct system of
embeddings.  Let $ T^*= \lng
A^*,\F^*,D^*,G^*,\phi^*\rng:=\limdir_I T_i$. Then $\F^*$ is
homogeneous, with $\phi^*$ testifying homogeneity, and the canonical
mapping $\pi_i:T_i\to T^*$ is a successful multi-embedding.}

\proof:   We first recall the definition of a direct limit.

An equivalence relation ${\equiv}$ is defined over $\bigdiscup_{i
\in I}\bar
T_i$ as follows: $a{\equiv} b\iff
(\exists
i\le j)(\pi_i^j(a)=b\vee \pi_i^j(b)=a)$. Conditions (i)--(ii) above
imply that $\equiv$ is indeed an equivalence relation. We define the
canonical map
$\pi_i(a)=[a]_{\equiv}$.  Next we set
$A^*=\bigdiscup_{i\in I}A_i/{\equiv}$ and observe the following:
\ppro \ctble Fact: {\sl For every infinite cardinal $\k$, if $I$ and
every $A_i$ are of cardinality $\le \k$, then $|A^*|\le\k$. }

We let $\F^*=\bigdiscup _{i\in I}\F_i/{\equiv}$,
$G^*=\bigdiscup _{i\in I}G_i/{\equiv}$ and $D^*=\bigdiscup_{i\in
I}D_i/{\equiv}$.

 For $x^*,y^*\in A^*, X^*\in \F^*, d^*\in D^*$ and $g^*\in G^*$ we
note:
\startitm
\itm  $x^*\in X^*$  iff there is some $i\in I$ and  $x\in A_i,X\in
\F_i$ such that $x\in X$ and $\pi_i(x)=x^*$, $\pi_i(X)=X^*$. 

\itm $g^*(x^*)=y^*$ iff $g(x)=y$ for some  $i\in I$ such that  $x\in
A_i,g\in G_i$  and $\pi_i(x)=x^*$, $\pi_i(y)=y^*$ and $\pi_i(g)=g^*$.
\ignore 
\itm $x\in \dom f^{d^*}, X^*\in \dom f^{d^*}, h^{d^*}(x^*)=y^*$ and
$f^{d^*}(X^*)=Y^*$ iff there is some $i\in I$ and $d\in D_i$ such
that
$\pi_i(d)=d^*$ and there are $x\in \dom h^d, X\in \dom f^d$, $y\in
A_i$ and $Y\in \F_i$ such that $\pi_i(x)=x^*$,
$\pi_i(X)=X^*,\pi_i(y)=y^*,\pi_i(Y)=Y^*$ and $h^d(x)=y$, $f^d(X)=Y$.

$d\in

\endignore 
\itm $\phi^*(d^*)=g^*$ iff there is $i\in I$ such that $\phi_i(d)=g$
and $\pi_i(d)=d^*, \pi_i(g)=g^*$.

We leave verification of this to  the  reader and that the following 

hold.

\startitm
\aitm $\F^*\su \Cal P(A^*)$ 

\aitm $G^*\su \aut(\F^*)$
\aitm $D^*=D(\F^*)$
\aitm $\phi^*:D^*\to G^*$ is a (total) satisfying homomorphism.
\aitm $\pi_j\pi_i^j=\pi_i$ for $i\le j$ in $I$

We  conclude that $\pi_i:T_i\to T^*$ is a successful embedding for
every
$i\in I$.

Homogeneity of $\F^*$ follows readily from
(c) and (d) above.\endproof\DirIsHom

\neusection
\beginsection{\S\number\secno\ Bi-universal homogeneous families}

The result proved in this section is the existence of a bi-universal
member in the class of homogeneous families over a given infinite  
set.

 Let us make the following definition:

\ppro \biembedding Definition:  {\sl We call an embedding of  
structures
$\pi:M\to N$ a {\it bi-embedding} if for every automorphism $g\in
\aut(M)$ there is an automorphism $g'\in \aut(N)$ such that
$\pi(g(x))=g'(\pi(x))$ for all $x\in M$.}

We observe that if $f:M\to N$ is a bi-embedding then $f$ induces an
embedding of $\aut(M)$ into the group of all restrictions to $f[M]$  
of
elements in the set-wise stabilizer of $f[M]$ in $\aut(N)$; that is,
an embedding as permutation groups (see [T1]). We can think of a
bi-embedding as a simultaneous embedding of both a structure and its
automorphism group.

\ppro\biuniversal Definition:  {\sl A structure $M^*$ in a class of
structures $K$ is {\it bi-universal} if for every structure $M\in K$
there is a bi-embedding $\pi:M\to M^*$.}

\ppro \remTruss Remarks:
\startitm
\itm The definition of embedding of permutation grpups (see [T1]) is
obtained by from this one by adding the condition that $\pi$ is onto.
\itm Example \example\ (1)
above indicates that if a bi-universal family $\F^*$ over a set $A^*$  
exists,
then for some $A\su A^*$ of cardinality $|A^*|$ the restrictions of
automorphisms of $\F^*$ to $A$ include the full symmetric group
$\Sym(A)$.

 \ppro \Hausdorff Lemma:  {\sl For every infinite $T=\lng
A,\F,D,G,\phi\rng$ there is a
set $B$ such  that  $|A|=|B|$ and a successful multi-embedding
$$\pi:T\to 
\lng A\discup B, \Cal P(A\discup B),D(A\discup B,\Cal P(A\discup
B)),\Sym(A\discup B),\phi'\rng$$
}

\proof We specify  the points of $B$. A point in $B$ is
a
finite function from the power set of a finite subset of $A$ to
$\{0,1\}$, namely $f\in B\iff f:\Cal P(D_f)\to \{0,1\}$ and $D_f\su
A$
is
finite. We let $\pi\rest A=\id$. For $X\in \F$ we define $\pi(X)$ as
follows: $\pi(X)=X\cup \{f\in
B: f(X\cap D_f)=1\}$. We let $\pi(\sigma)\rest A=\sigma$ and let
$\pi(\sigma)(f)=g\iff \sigma[D_f]=D_g\wedge f(X)=g(\sigma[X])$ for
all
$X\su D_f$. It is straightforward to verify that $\pi\rest\Sym(A)$
is a group monomorphism. We verify condition (c) in the definition of
successful embedding (definition \embedding\ above). Suppose  $X\su
A$
and $\sigma\in
\Sym(A)$ are given.

$\pi(\sigma)[\pi(X)]=$

$\sigma[X]\discup\pi(\sigma)[\{g\in B: g(X\cap D_g)=1\}]=$

$\sigma[X]\discup\{\pi(\sigma)(g):g\in B\wedge g(X\cap D_g)=1\}=$

$\sigma[X]\discup\{f\in B: f(\sigma[X]\cap D_f)=1\}=$

$\pi(\sigma[X])$

The definition of $\pi\rest D(\Cal P(A)$ is determined uniquely by
condition (b) in \embedding\ above. We need to specify $\phi'$ and
prove that (d) holds. For this we notice that:

\ppro \independence Claim: {\sl The family $\F=\{\pi(X):X\su A\}$ 

satisfies that  for every
finite
function $\tau:\F\to\{+,-\}$ the set $B_\tau=\bigcap_{X\in\tau^{-
1}(+)}\pi(X)\cap
\bigcap_{Y\in\tau^{-1}(-)}(A\discup B)\sm \pi(Y)$ has the same  
cardinality as
$A\discup B$.}

\proof The proof of this is well known.\endproof\independence

\ppro \moreIndependence Corollary: {\sl For every demand $d$ on $\F$  
there
is a permutation $\sigma\in\Sym(A\discup B)$ such that
$\sigma(x)=h^d(x)$ and $\pi(\sigma[X])=\pi[f^d(X)]$ for every $x\in
\dom h^d$ and $X\in\dom f^d$.}

\proof For every $\tau:\dom f^d\to \{+,-\}$ it holds that 

$|B_\tau|=|A\discup B|=|B'_\tau|$ where $B'_\tau=\bigcap_{X\in
\tau^{-1}(+)}\pi(f^d(X))\cap \bigcap_{X\in\tau^{-1}(-)}\pi(A\sm
f^d(X))$.
(This
means, informally, that every "cell" in the Venn
diagram of $\dom f^{\pi(d)}$ and every "cell" of the Venn diagram of
$\ran
f^{\pi(d)}$ is of cardinality $|A\discup B|$). Therefore it is
trivial to extend $h^d$
to a permutation that carries $B_\tau$ onto $B'_\tau$ for every
$\tau$.
\endproof{\moreIndependence}

Now let us define $\phi'(\pi(d))=\pi(\phi(d))$ for every $d\in \dom
\phi$ and for all $d\in D\sm\dom \phi$ let us pick by claims
\independence\ and \moreIndependence\ above a permutation
$\phi'(\pi(d))$ that extends
$\pi(d)$.\endproof\Hausdorff

\ppro \firstlim Theorem: {\sl Suppose $A_0$ is a given infinite set.
There
is a successful direct system of  embeddings ${\bf T}=\lng T_n: (n\in
\om); (\pi_m^n,\phi_m^n)\rng$ such that:
\startitm
\itm $A_n$ is of cardinality $|A_0|$
\itm $\F_n=\Cal P(A_n)$
\itm $G_n=\Sym(A_n)$. }

\proof  Let $T_0=\lng A_0,\Cal P(A_0),  
D(\F_0),\Sym(A_0),\phi_0:\{e\}\to
\{\id_{A_0}\}\rng$. Now use Lemma \Hausdorff\
inductively.\endproof\firstlim

\ppro \universal Theorem:  {\sl For every infinite set $A^*$ there is  
a
homogeneous family $\F^*\su\Cal P(A^*)$, and an infinite subset
$A\su A^*$ of cardinality $|A^*|$ such that $\Cal P(A)=\{X\cap
A:X\in \F^*\}$ and $\Sym(A)\su\{g\rest A:g\in \aut(\F^*)\}$.
Therefore {\it any} injection $f:A^*\to A$ induces a bi-embedding of
{\it every} family $\F\su \Cal P(A^*)$ (not necessarily homogeneous)
into $\Cal F^*$. In
particular, $\F^*$ is bi-universal in the class of all homogeneous
families over $A^*$. }

\proof By Theorem \firstlim\ there is a successful direct system of
embeddings ${\bf T}=\lng T_n: (n\in 
\om); (\pi_m^n,\phi_m^n)\rng$ such that:
\startitm
\itm $|A_n|=|A^*|$
\itm $\F_n=\Cal P(A_n)$
\itm $G_n=\Sym(A_n)$.

By Theorem \DirIsHom\ and the side remark \ctble\ it follows that the
family $\F^*$ obtained by the direct limit is a homogeneous family of
subsets of a set $A^{**}$ of size $|A^*|$, and we may assume that
$A^{**}=A^*$. The canonical map $\pi_0$ is a successful
multi-embedding, and therefore in particular a bi-embedding. Let $A$
be the image of $A_0$ under $\pi_0$. As $\F_0=\Cal P(A_0)$ and
$G_0=\Sym(A_0)$, we conclude that $\Cal P(A)=\{X\cap A:X\in \F^*\}$
and $\Sym(A)\su\{g\rest A:g\in \aut(\F^*)\}$. The Theorem is now
obvious.\endproof\universal
\ignore
One may wonder at this point if there is a homogeneous family $\F^*$
with the property that every other homogeneous family over a  
countable
set is isomorphic to one of its subfamilies. This would mean that
$\aut \F^*$ is universal in the class of all groups acting
homogeneously on some $\F\su\Cap P(A^*)$ with respect to the
definition in [T1].  We shall prove at the end
of the next section (theorem
\noUni\ below) that  no such family exists.
\endignore

\neusection
\beginsection{ \S\number\secno\ The number of isomorphism types of
homogeneous families over $\om$}

In this section we make a second use of the method of direct limits
as  introduced in Section 2 to determine the
number of isomorphism types of homogeneous families over a countable
set. It was
conjectured in [GGK] that this number is the maximal possible, namely
$2^\cont$. An isomorphism
between two families $\F_0\su \Cal P(A_0)$ and $\F_1\su \Cal P(A_1)$
is, of course, a 1-1 onto function $f:A_0\to A_1$ which satisfies
$X\in \F_0\iff f[X]\in \F_1$.

To obtain 

$2^\cont$ non isomorphic homogeneous families over a countable set,
it
is enough to obtain $2^\cont$ {\it different} such families; for then
dividing by isomorphism, the size of each class is $\cont$, and
therefore there are $\cont$ classes (see below).

The technique used to achieve this is embedding a family $\F\su \Cal
P(A)$ in a homogeneous family $\F'\su \Cal P(A^*)$ for $A^*\supseteq
A$ in such a way that $\{X\cap A: X\in \F'\}=\F$. In other words, we
will ``homogenize'' a family $\F$ ``without adding sets'' to $\F$.  
Thus,
starting with distinct $\F$-s we obtain distinct homogeneous $\F'$-s.

\ppro \homogenizing Lemma: {\sl There is a pair of countable sets
$A_0\su A^*$ (in fact, for {\it every} pair $A_0\su A^*$ of countable
sets satisfying $A^*\sm A_0$ infinite) such that for every family
$\F\su\Cal P(A_0)$ satisfying $\fin(A_0)\su F$ there is a homogeneous
family $\F'\su\Cal P(A^*)$ satisfying $\{X\cap A_0:X\in\F'\}=\F$}

This lemma determines the number of isomorphism types of homogeneous
families over a countable set:

\ppro \manyNonIso Corollary: {\sl There are $2^\cont$ isomorphism  
types of
homogeneous families over a countable set.}

\proof There are $2^\cont$ different families $\{\F_\a:\a<2^\cont\}$,
such that  $\fin (A_0)\su \F_\a\su\Cal P(A_0)$. For each $\F_\a$
there is, by \homogenizing, a homogeneous family $\F'_\a\su \Cal
P(A^*)$
that satisfies $\{X\cap A_0:X\in \F'_\a\}=\F_\a$. Therefore,
$\a\not=\b$ implies that $\F'_\a\not=\F'_\b$. Let us define an
equivalence relation over $2^\cont$: $\a\equiv\b\iff$ there is an
isomorphism between 

$\F'_\a$ and $\F'_\b$. There are at most $\cont$ many members in an
equivalence class $[\a]_\equiv$, as there are  $\cont$ many
permutations of $A^*$, and therefore at most $\cont$ many different
isomorphic
images of $\F'_\a$. As $\cont\times\cont=\cont$, while
$2^\cont>\cont$, there must be $2^\cont$ many equivalence classes
over $\equiv$, and therefore $2^\cont$ many isomorphism types of
homogeneous families over $A^*$. \endproof\manyNonIso

We prepare for the proof lemma \homogenizing. Before plunging into  
the
formalism, let us state the idea behind the proof. We use the set of
demands over a family and the free group associated with this set to
construct a successful extention in which the automorphisms act  
freely.
Thus, we can control sets in the orbit of an ``old'' set so that  
their
intersections with the ``old'' set is either finite or ``old''.

We need
some notation: Let $\FG(D)$ be the free group over the set $D=D(\F)$
for some family $\F$.
If $\F$ is countable, this group is also countable. We view $\FG(D)$
as the collection of all reduced words in the alphabet $C=D\cup  
\{d^{-
1}:d\in D\}$ (a word is reduced if there is no occurrence of $dd^{-
1}$ or $d^{-1}d$ in it) and the group operation, denoted by $\conc$,
is juxtaposition and cancellation (so $w_1\conc w_2$ is a reduced
word, and its length may be strictly smaller than $\lg w_1+\lg w_2$).
We let $c$ range over the alphabet $C$, and let $c^{-1}$ denote $d^{-
1}$
if $c=d$ or $d$ if $c=d^{-1}$. We denote by $e$ the unit of the free
group, which is the empty sequence $\emptyseq$. For convenient
discussion we also adopt the notation $h^c$ and $f^c$, by which we
mean $h^d$ and $f^d$ if $c=d$ and the respective inverses  
$(h^d)^{-1}$
and $(f^d)^{-1}$ otherwise. Now we can define:

\ppro \specialWord Definition:  {\sl Suppose that ${\bf T}=\lng
T_i:(i\in I);\pi_i^j\rng$ is a successful direct system of
multi-embeddings. For every $j\in I$:

 \startitm

\itm A homomorphism $\xi_j:\FG(D_i)\to G_j$ is defined by
$\xi_j(d):=\phi_j(\pi(d))$.

\itm We call a word $w=c_0\ldots c_k\in \FG(D_{j})$ {\it new} if  
$c_l$
is not in the range
of $\pi_i^{j}$ for all $l\le k$ and all $i<j$. A word $w\in\FG(D_j)$  
is {\it old} if
it is in the range of $\pi_i^j$ for some $i<j$.

\itm For a word $w\in\FG(D_j)$ and $X\in \F_i$ we define
what $f^w(X)$ is. Let $w=w_0w_1\ldots w_l$ where for each $k\le l$  
the
word $w_k$ is either new or old. For a new word
$w_k=c_k^0\ldots c_k^{l(k)}$ we denote by $f^{w_k}$ the composition
$f^{c_k^{l(k)}}\ldots f^{c_k^0}$. If this composition is empty, we  
say
that $f^{w_k}$ is not defined. If $w_k$ is old, then
$\xi_j(w)\in\aut(\F_j)$ and induces a 1-1 function $f^{w_k}:\F_j\to
\F_j$.
Let $f^w$ be the composition $f^{w_l}\ldots f^{w_0}$. If this
composition is empty, we say that $f^w$ is not defined.

\itm Analogously to the definition in $(3)$, we define $h^w$. }

\medskip
To prove lemma \homogenizing\ we need an expansion of the technique  
of
direct limits by some more structure. This is needed to enable us to
handle uncountably many demands by adding just countably many points.
We first define (a particular case of) inverse systems. Then we form
direct limits of inverse systems to obtain a pair of sets as required
by the lemma.

\ppro \invsys Definition:  {\sl a sequence ${\bf T}=\lng  
T^m:m<\om\rng$,
where $T^m=\lng A^m,\F^m,D^m,G^m,\phi^m\rng$, is called {\it an
inverse system} if:
\startitm
\itm $\F^m\su\Cal P(A^m)$, $D=D(A^m,\F^m)$, $G^m\su \aut (\F^m)$ and
$\phi^m:\FG(D^m)\to G^m$ is a partial satisfying homomorphism.
\itm $A^m$ and $\F^m$ are {\bf countable}

For $m\le m'$

\itm  $A^m\su A^{m'}$
\itm $\F^m\su\{X\cap A^m: X\in \F^{m'}\}$
\itm $G^m\su\{g\rest A^m:g\in G^{m'}, g\rest A^{m}\in \Sym(A^{m})\}$

For a demand $d\in D^{m'}$ we define $d\rest A^m$ iff $\dom h^d\cup
\ran
h^d\su A^m$ and for every distinct $X,Y\in\dom f^d\cup \ran f^d$ the
sets $X\cap A^m$ and $Y\cap A^m$ are distinct. When $d\rest A^m$ is
defined, $h^{d\rest A^m}=h^d$ and $f^{d\rest A^m}$ is obtained from
$f^d$ by replacing every $X\in \dom f^d\cup \ran f^d$ by $X\cap A^m$.
Clearly, when $d\rest A^m$ is defined, it belongs to $D^m$, and every
$d\in D^m$ equals $d'\rest A^m$ for some $d'\in D^{m'}$ by (3) and
(4).

If $w=c_0\ldots c_k\in FG(D^{m'})$ and $c_i\rest A^m$ is defined for
every $i\le k$, we define $w\rest A^m$ as $c_0\rest A^m\ldots  
c_k\rest
A^m$ (it is obvious what $c\rest A^m$ is). The restriction $\rest$ is
a partial homomorphism from $\FG(D^{m'})$ onto $\FG(D^m)$. The last
condition is

\itm If $d\in \dom \phi^{m'}$ and $d\rest A^m$ is defined, then  
$d\rest
A^m\in \dom \phi^m$ and $\phi^m(d\rest A^m)=\phi^{m'}(d)\rest A^m$  
(the
operation of $\phi^{m'}(d)$ on $A^m$ depends only on $d\rest A^m$  
when
$d\rest a^m$ is defined).}

\ppro \invLim Definition: {\sl Given an inverse system ${\bf T}=\lng
T^m:m<\om\rng$ we define the {\it inverse limit} $\liminv {\bf
T}=T^*=\lng A^*,\F^*,D^*,G^*,\phi^*\rng$ as follows:
\startitm
\aitm{} $A^*=\bigcup_m A^m$.

For every $x\in A^*$ let $m(x)$ be the  least $m$ such that $x\in  
A^{m}$.

\aitm{} $\F^*=\{X\su A^*: (X\cap A^m\in
\F^m)\hbox { for all but finitely many } m\}$.  For $X\in \F^*$ we  
let
$m(X)$ be the  least such that $X\cap A^m\in\F^m$  for every $m\ge  
m(X)$.

 We call $X\in\F^*$ {\it bounded} if $X\su A^m$ for some $m$.

\aitm{} $G^*=\{g\in\Sym(A^*):(g\rest
A^m\in G^m)\hbox { for all but finitely many } m\}$. Let $m(g)$ be  
the
least such that $g\rest
A^m\in G^m$ for every $m\ge m_g$.

 It is easy to verify that $G^*\su \aut (\F^*)$.

\aitm $D^*=D(\F^*)$

It is easy to verify that for every $d^*\in D^*$ there is some
$m(d^*)$ such that for all $m\ge m(d^*)$ it is true that $d^*\rest
A^m$ is defined, and $d^*\rest A^m\in D^m$.

\aitm $\phi^*(d^*)=\bigcup_{m\ge m_{d^*}} \phi^m(d^*\rest A^m)$ and  
is
defined iff $d^*\rest A^m\in \dom \phi^m$ for all $m\ge m_{d^*}$}

\ppro \multiInv Definition:  {\sl Suppose that ${\bf T}_0=\lng
T_0^m:m<\om\rng$ and ${\bf T}_1=\lng T_1^m:m<\om\rng$ are inverse  
systems,
and let $\liminv{\bf T}_0=T_0=\lng A_0,\F_0,D_0,G_0,\phi_1\rng$ and
$\liminv {\bf T}_1=T_1=\lng A_1,\F_1,D_1,G_1,\phi_1\rng$ be their
respective inverse limits. We call a sequence $\lng\m{m}\pi:T_0^m\to
T_1^m:m<\om\rng$ of multi-embeddings an {\it inverse system of
multi-embeddings} if for $m\le m'$ we have:
\startitm
\itm $\m{m'}\pi\rest A_0^m=\m{m}\pi\rest A^m_0$
\itm $\m{m'}\pi(X)\rest A^m_1=\m{m}\pi(X\cap A_0^m)$ for every $X\in
\F^{m'}_0$ for which $X\cap A_0^m\in\F_0^m$
\itm $\m{m'}\pi(g)\rest A^m_1=\m{m}\pi(g\rest A^m_0)$ for every $g\in
G_0^m$ for which $g\rest A_0^m\in G^m$



When  $\lng\m{m}\pi:T_0^m\to T_1^m:m<\om\rng$ is an inverse system
of multi-embeddings we   define a multi-embedding
$\pi=\liminv\m{m}\pi:{T}_0\to
{ T}_1$ as follows:
\item{} $\pi\rest A^*_0=\bigcup \m{m}\pi\rest A^m_0$
\item{} $\pi(X)=\bigcup _{m\ge m(X)}\m{m}\pi(X\cap A_0^m)$ for $X\in  
F^*_0$
\item{} $\pi(g)=\bigcup_{m\ge m_g} \m{m}\pi(g\cap A_0^m)$ for $g\in
G^*_0$

Call $\pi=\liminv \m{m}\pi$ a {\it multi-embedding of inverse
systems}.
 }
\ppro \invIsSucc Claim:  {\sl If  ${\bf T}_0=\lng
T_0^m:m<\om\}$ and  ${\bf T}_1=\lng T_1^m:m<\om\}$ are inverse
system and  $\lng\m{m}\pi:T_0^m\to T_1^m:m<\om\rng$ is  an inverse
system
of multi-embeddings such that every $\m{m}\pi$ is successful, then
$\pi=\liminv_m\m{m}\pi$ is also successful. }

\proof Suppose that $d\in D_0$ and we shall show that $\pi(d)\in\dom
\phi_1$. There is some $m_d$ such that for all $m\ge m_d$ the
restriction $d\rest A_m$ is defined. As $\m{m}\pi$ is successful,
$\m{m}\pi(d\rest A_m)$ belongs to $\dom \phi_1^m$ for $m\ge m_d$.
Therefore $\phi_1(\bigcup_{m\ge m_d}\m{m}\pi(d\rest
A_m)=\phi_1(\pi(d))$ exists and
belongs to $G_1$. \endproof\invIsSucc

\bigbreak\bigbreak

We shall construct a two dimensional system ${\bf T}=\lng
T_n^m:n,m<\om\rng$ and successful multi-embeddings
$\m{m}\pi_n^{n+1}:T_n^m\to T_{n+1}^m$ such that for every $n$,
\startitm
\itm  ${\bf T}_n=\lng T_n^m:m<\om\rng$ is an inverse system.
\itm  $\lng \m{m}\pi_n^{n+1}:m<\om\rng$
is an inverse system of successful  multi-embeddings.

Then a direct system will result: $T_n=\liminv {\bf T}_n$ and
$\pi_n^{n+1}=\liminv \m{m}\pi_n^{n+1}$.

\bigbreak
 Let $T_0^m=\lng m+1,\Cal P(m+1),D(\Cal
P(m+1)),\{\id\},\{(e,\id)\}\rng$.
  Clearly, $T_0=\liminv {\bf T}_0=\lng \om,\Cal
P(\om),D(\Cal P(\om)),\{\id\},\{(e,\id)\}\rng$.

Suppose now that $T_n=\liminv T_n^m$ is defined, where $T_n^m=\lng
A_n^m,\F_n^m,D_n^m,G_n^m,\phi_n^m\rng$, and that $\pi_{n-1}^n=\liminv
\m{m}\pi_{n-1}^n$ is also defined (when $n>0$)

We assume, for
simplicity, that $\pi_{n-1}^n\rest A_{n-1}=\id$ (if $n> 0$) and,
furthermore, identify $\FG(D_{n-1}^m)$ with its image under
$\m{m}\pi_{n-1}^n$, and write $\FG(D_{n-1}^m)\su\FG(D_n^m)$ as well
as $\FG(D_{n-1})\su\FG(D_n)$. Thus, the new words
of $\FG(D_n)$ coincide with $\FG(D_n\sm D_{n-1})$, and similarly for  
$\FG(D_n^m)$.

Let $\bigdis _m D_n^m$ be the disjoint union of $D_n^m$. 

We view $A_n$ as a subset of the following set $B_{n+1}= \{xw:x\in  
A_n,
w\in\FG(\bigdis_mD_n^m)\}$.  The expression $xw$ is the formal string
$xc_0\ldots c_w$ where
$w=c_0\ldots,c_k$, and $x$ is
identified with $xe$ (where $e$ is the empty string).


\ppro \anything Fact: {\sl $B_{n+1}$ is countable. }

The fact holds because  each $D_n^m$ is countable.

 Now define $B_{n+1}^m=\{xw:x\in A_n^m,w\in\FG(\bigdis _{m'\le
m}D_n^{m'})\}$. Clearly,
$A_n^m\su B_{n+1}^m$.

Next we define an operation $\xi_{n+1}(c):B_{n+1}\to B_{n+1}$ for
every $c\in D_n$ (there are, of course, uncountably many $c$-s!).

We want that $\xi_{n+1}(c)\rest B_{n+1}^n$ to depend only on $c\rest
A_n^n$ whenever $c\rest A_n^m$ is defined.

If $x\in\dom h^c$, we let $\xi_{n+1}(c)(x)=h^c(x)$.

For all other points
in $B_{n+1}$, we let
$\xi_{n+1}(c)(xw)=xw\conc (c\rest A_0^{m})$ if $m$ is the least such
that
$xw\in B_n^m$ {\bf and}  $c\rest A_0^m (\in C_1^m)$ is defined.

There is a unique extension of $\xi_{n+1}$ to a homomorphism from
$\FG(D_n)$ to $\Sym (B_{n+1})$, which we also call $\xi_{n+1}$.

\ppro \firstB Claim: {\sl For every $w\in\FG(D_n)$ there is some  
$m(w)$ such
that: 

\startitm 

\itm $B_{n+1}^m$ is invariant under $\xi_{n+1}(w)$ for
all $m\ge m_w$.  

\itm If $w\not=e$ then for every $xv\in B_{n+1}\sm
A_n^{m_w}$, we have $\xi_{n+1}(w)(xv)=xv\conc w\not=xv$.}

\proof (1) is clear from the definition. For (2) notice that if
$c\rest A_0^m$ is defined then the finitely many points in $\dom h^c$  
belong to
$A_0^m$. Then $\xi_{n+1}(c)(xv)=xv\conc (w\rest A_0^m)$.

From \firstB\ (2) it follows readily that $\xi_{n+1}$ is, in fact  a
{\it monomorphism}, as for every $w\in \FG(d_0)$ there is some $m_w$  
for
which $w\rest A_0^m$ is defined.

Let  $\m{m}\xi_{n+1}(c\rest A_n^m)=\xi_{n+1}(c)\rest B_{n+1}^n$ for  
all
$c\in D_n$ for which $c\rest A_n^m$ is defined.

Now we can define $A_{n+1}=\{\xi_{n+1}(w)(x):x\in A_n,  
w\in\FG(D_n)\}$
and $A_{n+1}^m=A_n\cap B_{n+1}^m=\{\m{m}\xi(w)(x):x\in
A_0^m,w\in\FG(D_0^m)$. (We remark that $A_{n+1}\not=B_{n+1}$, because  
when
$x\in\dom h^c$, the point $xc\notin A_{n+1}$).

Clearly, $A_{n+1}$ is invariant under  $\xi_{n+1}(w)$ for every $w\in
\FG(D_n)$, and also $A_{n+1}^{m_w}$ is, if $w\rest A_n^{m_w}$ is
defined.

Having defined $A_{n+1}$ we let $\m{m}\pi_n^{n+1}: A_n^m\to  
A_{n+1}^m$ be
the identity. Therefore also $\pi_n^{n+1}\rest A_n$ is the identity.

Now let us define $\m{m}\pi_n^{n+1}\rest \F_n^m$. For every
$X\in\F_n^m$ and $xw\in A_{n+1}^m$ we determine whether
$xw\in\m{m}\pi_n^{n+1}(X)$ by induction on the
length of $w$.

If $\lg w=0$ then necessarily $xw=x$, and we let
$x\in\m{m}\pi_0^1(X)\iff
x\in X$ for every $X\in \F_n^m$ and $x\in A_n^m$.

Suppose that this is done for all words of length $k$ and that $\lg  
wc=k+1$.

Distinguish two cases: when $c$ is old
and when $c$ is new.

First case: $c$ is old, namely $c\in C_{n-1}^{m}$ (this case does not
exist when $n=0$).  Here we have that
$\m{m}\xi_{n}(c)=\m{m}\phi_n(c)$ is
defined, and is an automorphism of $\F_n^m$. Let
$xwc\in\m{m}\pi_0^1(X)\iff xw\in\m{m}\pi_n^{n+1}(\m{m}\xi_{n-1}(c^{-
1})[X])$.

Second case: $c$ is new. Let $xwc\in\m{m}\pi_n^{n+1}(X)\iff
xw\in\m{m}\pi_n^{n+1}(f^{c{^-1}}(X))$.     In the right hand side we
mean that $f^{c^{-1}}(X)$ is defined and  
$xw\in\m{m}\pi_n^{n+1}(f^{c^{-
1}}(X))$.

Now we can set $\pi_n^{m+1}(X)=\bigcup_{m\ge  
m(X)}\m{m}\pi_n^{n+1}X\cap
A_n^m$.

\ppro \comm Fact: {\sl For every old $w\in \FG(D_n)$ and every $X\in  
F_n$
it holds that  
$\pi_n^{n+1}(\phi_n(w)[X])=\xi_{n+1}(w)[\pi_n^{n+1}(X)]$
(rule (c) in \embedding).}

The proof of the fact is straightforward using induction on word  
length.

\ppro \characterization Claim: {\sl For every $w\in\FG(D_n)$ and  
every
$X\in \F_n$ there is $m\ge m_w$ such that 

\startitm
\itm $f^w(X)$ is defined iff
$f^{w\rest A_n^m}(X\cap A_n^m)$ is defined
\itm  $f^{w^{-1}}(X)$ is
defined iff $f^{w^{-1}\rest A_n^m}(X\cap A_n^m)$ is defined 

\itm for every $x\in A_n$ with $m(x)\ge m$, $\xi_{n+1}(w)(x)=x(w\rest
A_n^{m(x)})\in\pi_n^{n+1}(X)\iff x\in f^{w^{-1}}(X)$ (where by $ x\in
f^{w^{-1}}(X)$ we mean that $f^{w^{-1}}(X)$ is defined {\bf and} $  
x\in
f^{w^{-1}}(X)$).
}

\proof If $f^w(X)$ is defined, then $f^{w\rest A_n^m}(X\cap A_n^m)$  
is
defined whenever $w\rest A_n^m$ is defined and equals $f^w(X)\cap
A_n^m$. Conversely, if $f^w(X)$ is not defined, then there is some
$m\ge m_w$ such that $X\cap A_n^m\not=Y\cap A_n^m$ for all $Y\in \dom
^w$ (if there is one $X$ for which $f^w(X)$ is not defined,
then  $\dom f^w$ is necessarily finite) and therefore $f^{w\rest
A_n^m}(X\cap A_n^m)$ is not defined.

From the definition of $\xi_{n+1}$ and  $m(x)\ge m_w$
 it follows that $\xi_{n+1}(x)=x(w\rest
A_n^m)$. From the definition of $\pi_n^{n+1}\rest \F_n$ it is
immediate that $x(w\rest A_n^m)\in\pi_n^{n+1}(X)\iff x\in
f^{w^{-1}}(X)$.\endproof\characterization

\ignore
As a corollary to the claim  we have:

If $f^w$ is not defined, then for some $m_0$ we have that for all  
$x\in A_n$
for which $m(x)\ge m_0$, $\xi_{n+1}(w)(x)\notin \pi_n^{n+1}(X)$ for  
any
$X\in \F_n$. Such a point of $A_{n+1}$ is called an {\it orphan}.
Trivially, if $xw$ is an orphan and $wv$ is an extension of $w$, then
 $xwv$ is also an orphan. (One remains an orphan even when one grows
longer).
\endignore

\ppro \secondB Fact: {\sl $\pi_n^{n+1}(X)\cap A_{n+1}^m$ depends only  
on
$X\cap A_n^m$ whenever $X\cap A_n^m\in\F_n^m$.}\endproof\secondB

Now we can define  

$\F_{n+1}=\{\xi_{n+1}(w)[\pi_n^{n+1}(X)]:X\in\F_n,w\in\FG(D_n)\}$.

Let  

$\F_{n+1}^m=\{\m{m}\xi_{n+1}(w)(X):X\in\F_{n+1}^m,w\in\FG(D_n^m)\}$.

\ppro \therdB Claim: {\sl $\F_{n+1}^m$ is countable for every $m$. 

\proof The fact follows by the countability of $\FG(D_n^m)$ and 

\secondB.} \endproof\therdB

We finished defining ${\bf T}_{n+1}$ and $\lng
\m{m}\pi_n^{n+1}:m<\om\rng$, and verified that ${\bf T}_{n+1}$ is an
inverse system, that  $\lng
\m{m}\pi_n^{n+1}:m<\om\rng$ is an inverse system of successful
multi-embedding and that, consequently, $\pi_n^{n+1}:T_n\to T_{n+1}$
is a multi-embedding of inverse systems.

\ignore
\ppro \theStar Claim:  {\sl Suppose that $e\not=w\in FG(D_n\sm  
D_{n-1})$ is
new
and $X\in\F_n$. If $f^w(X)$ is defined, then
$\xi_{n+1}(w)[\pi_n^{n+1}(X)]=\pi_n^{n+1}(f^w(X))$. If $f^w(X)$ is  
not
defined, then $\xi_{n+1}(w)[\pi_n^{n+1}(X)]\cap A_n$ is bounded [and
$\xi_{n+1}(w)[\pi_n^{n+1}(X)]\not=\pi_n^{n+1}(Y)$ for all $Y\in  
\F_n$.]}

\proof We prove the first part of the claim by induction on $\lg w$.
Suppose that $wc$ is of length $k+1$.

Suppose that $f^{wc}(X)$ is defined. Then
also $f^w(X)=Y$ is defined, unless when $w=e$, a case in which we let
$Y=X$. So $f^{wc}(X)=f^c(Y)$. By the induction hypothesis we know
that $xv\in\pi_n^{n+1}(X)\iff \xi_{n+1}(w)(xv)\in\pi_n^{n+1}(X)$ for
every $xv\in A_{n+1}$. We show that $xv\in\pi_n^{n+1}(Y)\iff
\xi_{n+1}(c)(xv)\in\pi_n^{n+1}$ for every $xv\in A_{n+1}$. This holds
when $xv=x$ and $x\in \dom h^c$ by the definition of $\xi_{n+1}(c)$.
For $xv\notin \dom h^c$  $m$ be the first such that $c\rest A_n^m$
is defined, $xv\in A_{n+1}^m$  and $Y\cap A_n^m\in\F_n^m$. So
$\xi_{n+1}(c)(xv)=xv\conc c\rest A_n^m$. Now look at the definition  
of
$\pi_n^{n+1}(f^c(Y))$.

For the second part we do not need induction.  Suppose that
$f^{w}(X)$ is not defined.  Can it happen that
$f^{w\rest A_n^m}(X\cap A_n^m)$ is defined nevertheless? The answer  
is
yes, but only boundedly often.  There is a large enough $m_0$ such  
that
$X\cap A_n^m\not=Y\cap A_n^m$ for all $Y\in \dom f^w$. Then from
$m_0$ onwards, $f^{w\rest A_n^m}(X\cap A_n^m)$ is not defined.

Let $m\ge m_0$ be large enough so that $w\rest A_n^m$ is defined.  
Then
$A_{n+1}^m$ is invariant under $\xi_{n+1}(w)$.  We show that for  
every
$x\in A_n\sm A_n^m$, $\xi_{n+1}(w^{-1})(x)\notin \pi_n^{n+1}(X)$.  
This
is enough to show that $\xi_{n+1}(w)[\pi_n^{n+1}(X)]\cap A_n$ is
bounded.  Let $x\in A_n\sm A_0^m$, and let $m_1$ be the first such
that $x\in A_n^{m_1}$. As $w\rest A_0^m$ is defined, also  
$w^{-1}\rest
A_n^{m_1}$ is defined, and therefore
$\xi_{n+1}(w^{-1})(x)=x(w^{-1}\rest A_n^{m_1})$.
This cannot belong to $\pi_n^{n+1}(X)$ because $X\cap A_n^{m_1}\notin
\ran
f^{w^{-1}\rest A_n^{m_1}}$.
\endignore

 Let $T^*=\lng
A^*,\F^*,D^*,G^*,\phi^*\rng=\limdir T_n$.
We show that the conclusion of lemma \homogenizing\ holds for the  
pair
of sets $A_0$ and $ A^*$. Clearly, these sets are countable and  
$A_0\su
A^*$. So all we need is:

\ppro \final Claim: {\sl For every family $\F\su\Cal P(A_0)$ which  
includes
$\fin(A_0)$ there is a homogeneous family $\F'\su \F^*$ such that
$\F'\rest A_0=\F$.}

\proof  Suppose that $\F\su
\Cal P(A_0)$ is a family which includes $\fin(A_0)$. We work by
induction on $n$ and define $\F'_n\su\Cal P(A_n)$ for every $n$:
\startitm
\itm $\F'_0=\F$. 

\itm $\F'_{n+1}=\{\xi_n(w)[\pi_n^{n+1}(X)]:w\in\FG(D(\F'_n)),X\in
\F'_n\}$

\medbreak

Let $\F'=\{\pi_n(X):X\in\F'_n\}$.

We claim that
\startitm
\aitm $\F'\su \F^*$ and  $\phi^*\rest D(\F')$ testifies that $\F'$ is
homogeneous.
\aitm $\{X\cap A_0:X\in\F'\}=\F$.

To prove (a) suppose that $d\in D(\F')$ is a demand. Then there is  
some
$n$ and a demand $d_n\in D(\F'_n)$ such that $\pi_n(d_n)=d$. As
$\pi_n^{n+1}$ is successful,  
$\xi_n(d)=\phi_{n+1}(\pi_n^{n+1}(d_n))=:g$
is 

defined. Now $\pi_{n+1}(g)=\phi^*(d)$  satisfies $d$ and is an
automorphism of $\F^*$. Why is it also an automorphism of $\F'$?
Because of (2) above.

To prove (b) we notice that it is enough to prove by induction that  
for
every $n$ and $X\in\F'_{n+1}$, we have
\item{$(*)_n$}  $X\cap A_n\in\F'_n$ or is
bounded.

 For then it follows by induction that $X\cap A_0\in \F$ for
every $n$ and $X\in \F_n$: if $X\cap A_n\in\F_n$ we have that $X\cap
A_0\in\F$ by the induction; if $X\cap A_n$ is bounded, then $X\cap
A_0$ is finite and again in $\F$.

So let us prove $(*)_n$.
 We have to show that for every $w\in \FG(D_n)$
and every $X\in \F_n$ the set $\xi_{n+1}[\pi_n^{n+1}(X)]\cap A_n$
belongs to $\F'_n$ or is bounded.
 We show something stronger.

\ignore
First let us extend the definition of $f^w$ from new $w$ to all
$w\in\FG(D_n)$.  Suppose that $w=c_0c_1\ldots c_k\in FG(D_n)$. Let
rewrite $w=w_0w_1\ldots w_l$ where for each $i\le l$ the word $w_i$  
is
either new (belongs to $\FG(D_n\sm D_{n-1})$ or old (belongs to
i$\FG(D_{n-1})$). For every  old $w_i$ we have that $\phi_n(w)\in\aut
\F`_n$. Let $f^{w_i}:\F'_n\to\F'_n$ be the permutation determined by
$w_i$. For every new $w_j$ let   $f^{w_j}$ be as before. Finally, let
$f^w$ be the composition $f^{w_l}f^{w_{l-1}}\ldots f^{w_0}$.
\endignore

\item{$(**)_n$} For every $X\in \F'_n$ and $w\in\FG(D(\F'_n))$ if
$f^w(X)$ is defined then
$\xi_{n+1}[\pi_n^{n+1}(X)]=\pi_n^{n+1}(f^w(X))$ (and therefore
$\xi_{n+1}[\pi_n^{n+1}(X)]\cap A_n= \pi_n^{n+1}(f^w(X))\cap
A_n=f^w(X)\in \F'_n$). If $f^w(X)$ is not defined, then
$\xi_{n+1}[\pi_n^{n+1}(X)]\cap A_n$ is bounded.

Suppose first that $f^w(X)$ is defined. Then obviously it belongs to
$\F'_n$, because $w\in\FG(D(\F'_n))$. It is easy to check that
$\xi_{n+1}(w)(xv)\in\pi_n^{n+1}(f^w(X))\iff xv\in\pi_n^{n+1}(X)$.

So assume that $f^w(X)$ is not defined, and we want to prove that
$\xi_{n+1}(w)[\pi_n^{n+1}(X)]\cap A_n$ is bounded.

If  $f^w(X)$ is not defined, then $X\notin\ran f^{w^{-1}}$. It is
sufficient to see that the set

$$\{x\in A_n:\xi_{n+1}(w^{-1})(x)\in\pi_n^{n+1}(X)\}$$

is bounded, because this set equals $\xi_{n+1}(w)[\pi_n^{n+1}(X)]\cap
A_n$. By \characterization\ there is a large enough $m> m(w)$ such
that for all $x\in A_n$ with $m(x)\ge m$ we have that

$$\xi_{n+1}(w)(x)=\m{m(x)}\xi_{n+1}(x)=x(w\rest A_n^{m(x)})\in
\m{m(x)}\pi_n^{n+1}(X\cap A_n^{m(x)})\iff x\in f^{w^{-1}}(X)$$

But $f^{w^{-1}}(X)$ is not defined, and therefore
$\xi_{n+1}(w)(x)\notin\pi_n^{n+1}(X)$ for all $x\in A_n$ with $m(x)>
m$, which is what we wanted. \endproof\homogenizing

We give a corollary of this proof.

\ppro \freeGroups Corollary: {\sl There is a collection of $2^\cont$
permutation groups over $\om$, $\lng G_\a:\a<2^\cont\rng$ such that:
\startitm
\itm Every $G_\a$ is isomorphic to the free group on $\cont$
generators.
\itm Every $G_\a$ testifies the homogeneity of some family
$\F_\a\su\Cal P(\om)$
\itm If $\a<\b<\cont$, then  $G_\a$ and
$G_\b$ are not isomorphic as permuatatio groups.}

\proof We have shown that there are $2^\cont$ many homogeneous
sub-families of $F^*$, $F'_\a$ for $\a<2^\cont$. The restriction of
$\phi^*$ to $\FG(D(A^*,\F_\a))$ is a monomorphism of the free group
over a set of cardinality $\cont$ into $G^*$ which testifies
homogeneity of $\F_\a$.  This gives us $2^\cont$ {\it different}
groups satisfying (1) and (2) in the corollary. To obtain (3), divide
by the relation ``isomorphic via a permutation of $\om$'', and pick a
member from every equivalence class.  As in each class there are
$\cont$ many members at the most, we get that there are $2^\cont$
classes.\endproof\freeGroups

\medbreak

We now wish to show that there is no homogeneous family over $\om$
such that every homogeneous family over $\om$ is isomorphic to one of
its subfamilies. This will follow from the next lemma about the  
number
of pairwise incompatible homogeneous families over a countable set.
Two  families over $\om$ are {\it incompatible} if for some
$X\su \om$ the set $X$ belongs to one family while the set $\om\sm X$
belongs to the other. For every $X\su \om$ let us denote $X^0:=X$ and
$X^1:=\om\sm X$.

\ppro \pairwise Lemma: {\sl  There is a collection  
$\{F_\a:\a<2^\cont\}$ of
pairwise incompatible homogeneous families over $\om$.}

\ppro \noUni Corollary: There is no homogeneous family over $\om$  
such
that every homogeneous family over $\om$ is isomorphic to one of its
subfamilies.

\proof\ (of Corollary) Suppose to the contrary that $\F^*$ is a
homogeneous family over $\om$ with this property. By Lemma \pairwise\
pick a collection $\{F_\a:\a<2^\cont\}$ of pairwise incompatible
homogeneous families over $\om$. For each $\a< 2^\cont$ fix a
permutation $\sigma_\a$ which embeds $\F_\a$ in $\F^*$. By the pigeon
hole principle there are $\a<\b<2^\cont$ and a permutation $\sigma$
such that $\sigma_\a=\sigma_\b=\sigma$. As $\F_\a$ and $\F_\b$ are
incompatible, let us find a set $X\su\om$ such that $X^0\in F_\a$ and
$ X^1\in \F_\b$.  Now $\sigma_\a(X^0)=\sigma(X^0)\in \F^*$, and
$\sigma_\b(X^1)=\sigma(X^1)\in \F^*$. This means that in $\F^*$ there
is a set and its complement. This contradicts Theorem  
\classification\
that states that there is no homogeneous family over $\om$ that
contains a set and its complement.\endproof\noUni

We prove now lemma \pairwise.

\proof We use the direct system of inverse systems from the proof of
lemma \homogenizing. The pairwise disjoint families will be over  
$A^*$
rather than over $\om$, but as this is a countable set this makes no
difference.

Let the variable $\eta$ range over the set of all functions  
$\eta:\Cal
P(A_0)\to 2$ which satisfy $\eta(X)+\eta(A_0\sm X)=1$ for all $X\su
A_0$.  These are functions that select exactly one element from each
pair of a set and its complement (for example, characteristic
functions of ultra filters). There are $2^\cont$ such functions.

For every function $\eta:\Cal P(A_0)\to 2$ as above let
$\F^0_\eta=\{X\su A_0:\eta(X)=1\}$. The collection
$\{\F_\eta:\eta:\Cal P(A_0)\to 2\}$ is a collection of $2^\cont$
pairwise incompatible families over $A_0$. For every $m<\om$ let 

 $\F_\eta^{0,m}$  be the projection of  $\F_\eta$ on  $A_0^m$.

We know that for every $\F_\eta$ there is a homogeneous family
$\F'\eta$ over $A^*$ whose projection on $A_0$  equals $\F_\eta$
(modulo finite sets). However, it is NOT true that  
$\{\F_\eta':\;\,\eta:\Cal P(A_0)\to 2\}$ is a collection of pairwise  
incompatible  families. In fact,  $\pi_0^1(X^0)\cap \pi_0^1(X^1)$ is  
not empty for every  $X\su A_0$.

What we shall do now is refine the extension operation is such a way
that not only the projection on $A_0$ is preserved, but also the
disjointness of $X^0$ and $X^1$. This will be achieved by removing
some of the points of $A^*$.

We define by induction on $n$ a subset $\bar D_n\su D_n$ and a subset
$E_n\su A_n$. Restricting ourselves to the points of $E=\bigcup_nE_n$
will provide the desired conservation property.

Let $E_0=A_0$. Let $\bar D_0=\bigcup_\eta D(E_0,\F_\eta^0)$.

\ppro \breath Fact:{\sl  If  $d\in \bar D_0$ then for no  $X\su A_0$  
is it
true that both $X^0,X^1$ belong to $\ran f^d$.}

We remove, thus, from the collection of demands all demands which
mention simultaneously a set and its complement in their range.

Let us now define $E_1$ as follows: 

 $$E_{1}= \{xw:x\in E_0,
w=c_0\ldots c_k\in\FG(\bigdis_m\bar D_0^m)\ \&\  x\notin\dom  
f^{c_0}\}$$

The variation on to the proof of \homogenizing\ is that only
a proper subset of words is being used. Hence,  $E_1\su A_1$

\ppro \disjoint Claim: {\sl For every  $X\su A_0$ it holds that
$\pi_0^1(X^0)\cap \pi_0^1(X^1)\cap E_1=\emptyset$.}

\proof By induction on the length of  $w\in FG(\bigdiscup_m \bar
D_0^m)$ we shall see that  $xw\notin \pi_0^1(X^0)\cap \pi_0^1(X^1)$.

If $\lg w=0$ then $xw=x\in E_0=A_0$. As $\pi_0^1(X)\cap A_0=X$ for  
all
$X$, it follows that  $x\notin \pi_0^1(X^0)\cap \pi_0^1(X^1)$.

Now suppose that $\lg wc=k+1$. By the definition of the $\in$  
relation
over the set $A_1$ we know that $xwc\in \pi_0^1(X^0)$ iff there is
some $Y$ such that $xw\in\pi_0^1(Y)$ and $f^c(Y)=X^0$. Similarly,
$xwc\in \pi_0^1(X^1)$ iff there is some $Z$ such that
$xw\in\pi_0^1(Z)$ and $f^c(Z)=X^1$. But $X^0$ and $X^1$ cannot both
appear in $\ran f^c$ because $c\in \bar D_0^m$. Therefore $xwc$ is  
not
in the intersection. \endproof\disjoint

Now we should notice that $E_1$ is invariant under $\xi_1(w)$ for all
$w\in\FG(\bar D_0)$. Also, for every $w\in\FG(\bar D_0)$ and every
$X\su E_0$ it holds that $\xi_1(w)[\pi_0^1(X^0)]\cap
\xi_1(w)[\pi_0^1(X^1)\cap E_1=\emptyset$.

Let  $\bar \F_1=\{\xi_1(w)[\pi_0^1(X)]:X\in\bar F_0,w\in\bar D_0\}$.

We proceed by induction on $n$, defining $\bar D_n$ and $E_{n+1}$ for
all $n>0$.

First, let us view each $\eta:\Cal P(E_0)\to 2$ as a partial function
$\eta:\bar\F_1\to 2$ by replacing every $X\su E_0$ by $\pi_0^1(X)$.
Next extend each $\eta$ to contain $\bar \F_1$ in its domain,
demanding that 

$$\eta (\xi_1(w)[\pi_0^1(X)])=\eta(X)$$

We refer to the resulting extended function also as $\eta$ to avoid
cumbersome notation. For every $\eta$ let $\bar \F_{\eta,1}=\{X\in
\bar \F_1:\eta(X)=1\}$.

Now define  $\bar D_1=\bigcup_\eta D(\bar F_{\eta,1})$.

Define  $E_{n+1}$ e`z $\bar \F_{n+1}$ as before. We should
check the following:

\ppro \disjointagain Claim: {\sl For all  $X\in\bar \F_n$ it holds  
that  $\pi_n^{n+1}(X^0)\cap
\pi_n^{n+1}(X^1)\cap E_{n+1}=\emptyset$.}

\proof By induction of word length. The case which should be added to
the proof of \disjoint\ is the case when $c$ as old, and is easily  
verified.

Having done the induction, we set $E=\bigcup _n E_n$. For every
$\eta:\Cal P(E_0)\to 2$ let $\F'_\eta$ be the homogeneous family
obtained from $\F_\eta$ as in the proof of \homogenizing. The reader
will verify that 

\startitm
\itm For every  $X\su E_0$ it holds that
$\pi_0(X^0)\cap\pi_0(X^1)\cap E=\emptyset$
\itm For every  $\eta:\Cal P(E_0)\to 2$ the family  $\F'_\eta\rest E$
is homogeneous.

This completes the proof. \endproof\pairwise

\bigbreak
\bigbreak
\noindent
{\bf References}

[GGK] M.~Goldstern, R.~Grossberg and M.~Kojman, {\sl Infinite
homogeneous bipartite graphs with unequal sides}, Discrete.

[KjSh1] M. Kojman and S. Shelah, {\sl 

Non existence of universal  oreders in many cardinals}, {\bf
Journal of Symbolic Logic} 57 (1992) 875--891.

[KjSh2]  M. Kojman and S. Shelah, {\sl The universality
spectrum of stable unsuperstable theories},  {\bf Annals of
Pure and Applied Logic}, 58 (1992) 57--72.

[KjSh3], M.Kojman and S. Shelah,  {\sl Universal Abelian Groups},
 Israel Journal of Math, to appear

[LW] A.~H.~Lachlan and R~.E.~Woodrow,
	{\sl  Countable Ultrahomogeneous Undirected Graphs}, 

Trans. Amer.	Math. Soc. 262 (1980) 51--94.

[MSST] A.~Mekler, R.Schipperus, S.~Shelah and J.~K.~Truss {The random
graph and automorphisms of the rational world} Bull. London Math.  
Soc.
25 (1993) 343-346

[Sh-c] S.~Shelah, {\bf Classification theory: and the number of non-
isomorphic models}, revised, {\it North Holland Publ. Co.}, Studies  
in
Logic and the Foundation of Math vol. 92, 1990, 705 + xxxiv.

[Sh-266] S.~Shelah, {\sl Borel Sets with large Squares}, in  
preparation.

[T1] J.~K.~Truss, {\sl Embeddings of Infinite Permutation Groups} in
{\it Proceedings of Groups --- St Andrews 1985} London Math.
Soc. Lecture Note Series no. 121 (Cambridge University Press 1986).
pp. 355--351

[T2] J.~K.~Truss 

{\sl The group  of the countable universal graph}  Math. Proc. Camb.
Phil. Soc (1985) {\bf 98}, 213--245

[T3] J.~K.~Truss 

{\sl The group of almost automorphisms of the
countable universal graph} Math. Proc. Camb. Phil. Soc (1989) {\bf
105}, 223--236

\bye